\documentclass[11pt,a4paper,twoside]{amsart}
\usepackage{amssymb}
\usepackage{mathtools}
\usepackage{graphics}
\usepackage{bm}
\usepackage{tikz}
\usepackage{xstring}
\usepackage{wasysym}
\newtheorem{theo}{\bf Theorem}

\theoremstyle{remark}

\newcommand{\mbb}{\mathbb}

\newcommand{\la}{\leftarrow}

\DeclareMathOperator{\maj}{\tt maj}
\DeclareMathOperator{\amaj}{\tt amaj}
\DeclareMathOperator{\des}{\tt des}
\DeclareMathOperator{\ides}{\tt ides}
\DeclareMathOperator{\asc}{\tt asc}
\DeclareMathOperator{\row}{\tt row}

\DeclareMathOperator{\DES}{\tt DES}
\DeclareMathOperator{\IDES}{\tt IDES}
\DeclareMathOperator{\ASC}{\tt ASC}
\DeclareMathOperator{\ROW}{\tt ROW}

\newcommand{\exA}{
	\draw[ultra thick] (0.5,0.5)--(0.5,6.5);
	\draw[thick] (0.5,0.5)--(6.5,0.5);
	\draw[thick] (6.5,0.5)--(6.5,6.5);
	\draw[thick] (0.5,6.5)--(6.5,6.5);
	\node[circle,fill] at (1,3){};
	\node[circle,fill] at (2,2){};
	\node[circle,fill] at (3,5){};
	\node[circle,fill] at (4,6){};
	\node[circle,fill] at (5,4){};
	\node[circle,fill] at (6,1){};
}
\newcommand{\exC} {
	\node at (7,5.5){$\la 1$};
	\node at (7,2.5){$\la 2$};
	\node at (8,3){};
}

\newcommand{\exI} {
	\draw[thick] (0.5,0.5)--(0.5,9.5);
	\draw[thick] (0.5,0.5)--(9.5,0.5);
	\draw[thick] (9.5,0.5)--(9.5,9.5);
	\draw[thick] (0.5,9.5)--(9.5,9.5);
	\node[circle,fill] at (1,5){};
	\node[circle,fill] at (2,2){};
	\node[circle,fill] at (3,7){};
	\node[circle,fill] at (4,9){};
	\node[circle,fill] at (5,6){};
	\node[circle,fill] at (6,1){};
	\node[circle,fill] at (7,8){};
	\node[circle,fill] at (8,4){};
	\node[circle,fill] at (9,3){};
}

\newcommand{\exF} {
	\node[circle,fill] at (1,1){};
	\node[circle,fill] at (2,2){};
	\node[circle,fill] at (3,3){};
	\node[circle,fill] at (4,1){};
	\node[circle,fill] at (5,3){};
	\node[circle,fill] at (6,5){};

	\draw[thick] (0.5,0.5) to (6.5, 0.5);
	\draw[thick] (6.5,0.5) to (6.5, 6.5);

	\draw[thick] (0.5,0.5) to (0.5, 1.5);
	\draw[thick] (1.5,1.5) to (1.5, 2.5);
	\draw[thick] (2.5,2.5) to (2.5, 3.5);
	\draw[thick] (3.5,3.5) to (3.5, 4.5);
	\draw[thick] (4.5,4.5) to (4.5, 5.5);
	\draw[thick] (5.5,5.5) to (5.5, 6.5);

	\draw[thick] (0.5,1.5) to (1.5, 1.5);
	\draw[thick] (1.5,2.5) to (2.5, 2.5);
	\draw[thick] (2.5,3.5) to (3.5, 3.5);
	\draw[thick] (3.5,4.5) to (4.5, 4.5);
	\draw[thick] (4.5,5.5) to (5.5, 5.5);
	\draw[thick] (5.5,6.5) to (6.5, 6.5);
}

\newcommand{\exH} {
	\draw[dashed] (6.5,6.5) to (6.5,7.5);
	\draw[dashed] (7.5,7.5) to (7.5,8.5);
	\draw[dashed] (8.5,8.5) to (8.5,9.5);
	\draw[dashed] (6.5,7.5) to (7.5,7.5);
	\draw[dashed] (7.5,8.5) to (8.5,8.5);
	\draw[dashed] (8.5,9.5) to (9.5,9.5);

	\draw[dashed] (6.5,0.5) to (9.5,0.5);
	\draw[dashed] (9.5,9.5) to (9.5,0.5);

	\node at (10, 1){$\la$};
	\node at (10, 4){$\la$};
	\node at (10, 7){$\la$};
}
\newcommand{\exJ} {
	\draw[thick] (6.5,6.5) to (6.5,7.5);
	\draw[thick] (7.5,7.5) to (7.5,8.5);
	\draw[thick] (8.5,8.5) to (8.5,9.5);
	\draw[thick] (6.5,7.5) to (7.5,7.5);
	\draw[thick] (7.5,8.5) to (8.5,8.5);
	\draw[thick] (8.5,9.5) to (9.5,9.5);

	\draw[thick] (6.5,0.5) to (9.5,0.5);
	\draw[thick] (9.5,9.5) to (9.5,0.5);

	\node[circle,fill] at (7, 1){};
	\node[circle,fill] at (8, 4){};
	\node[circle,fill] at (9, 7){};
	\draw[white,ultra thick] (6.5,0.5) to (6.5,6.5);
}

\begin{document}
\title{The double Eulerian polynomial and inversion tables}
\author{Erik Aas}
\address{Department of Mathematics, Royal Institute of Technology \\
  SE-100 44 Stockholm, Sweden}
\email{eaas@kth.se}
\date{January 2014}

\maketitle
{{\bf Abstract.} We show that the pair $(\des, \ides)$ of statistics on the set of permutations has the same distribution as the pair $(\asc, \row)$ of statistics on the set of inversion tables, proving a conjecture of Visontai. The common generating function of these pairs is the {\it double Eulerian polynomial}.}
\vspace{1cm}

\subsection*{The double Eulerian polynomial}
The double Eulerian polynomial $A_n(u,v)$ enumerates the number of descents of a permutation and its inverse,
\[
A_n(t, s) = \sum_{\pi\in\mathbb{S}_n} u^{\des(\pi)}v^{\des(\pi^{-1})}.
\]

It is a natural generalization of the classical Eulerian polynomial $A_n(u, 1)$. The latter polynomial is well-known to be positive in the basis $(u^i(1+u)^{n-i})_{i=0} ^n$. This has been proved in several ways; notably by geometric means \cite{BR}, and by an elegant bijective argument by Foata and Strehl (see \cite{P} for an excellent exposition).

There is no analogous result for the double Eulerian polynomial, though there is a conjectured one by Gessel \cite{B}: $A_n(u,v)$ should be integral and positive in the basis $((uv)^i(u+v)^j(1+uv)^{n-2i-j})_{i,j}$. Visontai \cite{V} gave explicit formulas for the coordinates of $A_n(u,v)$ in this basis, but was unable to prove that they are positive, nor that they are integers. He also conjectured a new way of defining $A_n(u,v)$, as follows (this is our main theorem).

\begin{theo}
\label{th_main}
	For all $n$, 
		\[
			A_n(u,v) = \sum_{e \in \mathbb{I}_n} u^{\asc(e)} v^{\row(e)}.
		\]
\end{theo}

Here, $\mbb{I}_n$ is the set of inversion tables of length $n$, and $(\asc, \row)$ are two statistics, all defined below. We will spend the remainder of this note proving Theorem \ref{th_main}, after giving the necessary definitions.

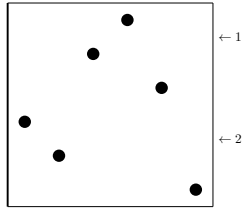
\begin{figure}
\scalebox{0.45}{
\begin{tikzpicture}
	\exA
	\exC
\end{tikzpicture}
}
\caption{A permutation $\pi' = 325641$ together with an $6+1$-tuple $(x_1, \dots, x_7) = (0,0,2,0,0,1,0)$ determines the child $\pi$ of $\pi'$ given in Figure \ref{fi_pi}. For all the examples in Figures 1-4, the fixed set $S$ can be taken to be $\{1,4, 7,8 \}$.}
\label{fi_pip}
\end{figure}

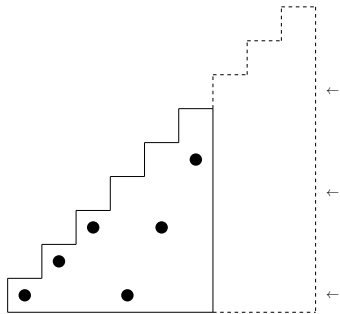
\begin{figure}
\scalebox{0.45} {
\begin{tikzpicture}
\exF
\exH
\end{tikzpicture}
}
\caption{An inversion table $e'=123135$ and a set $T = \{1,4,7\}$ determines the child $e$ of $e'$ given in Figure \ref{fi_e}.}
\label{fi_ep}
\end{figure}

\begin{figure}
\scalebox{0.45} {
\begin{tikzpicture}
\exF
\exJ
\end{tikzpicture}
}
\caption{The child $e=123135147$ to $e'$ referred to in Figure \ref{fi_ep}.}
\label{fi_e}
\end{figure}
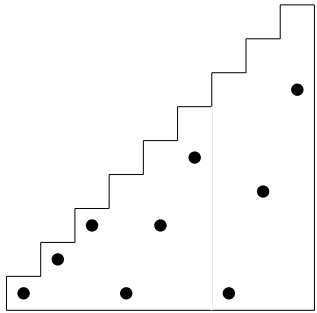

We identify permutations $w$ of length $n$ with words $w_1\dots w_n$, and with permutation diagrams $\{(i,w_i) : 1 \leq i \leq n\}$, which we read as Cartesian coordinates: a point $(x, y)$ refers to a point $x$ steps to the right and $y$ steps up from $(0,0)$. See Figures \ref{fi_pip} and \ref{fi_pi} for examples. An inversion table $e = e_1\dots e_n$ of length $n$ is any sequence of positive integers satisfying $1 \leq e_i \leq i$ for all $i$ (this differs slightly from the notation in \cite{S}, which we otherwise follow). We will identify these with marked staircases, examples of which are given in Figures \ref{fi_ep} and \ref{fi_e}. For a permutation $w$, let $\DES(w) = \{i : w_i > w_{i+1}\}$ be the {\it set} of descent positions, and $\IDES(w) = \DES(\pi^{-1})$ be the descent set of the inverse permutation. For an inversion table $e$, we define $\ASC(e) = \{i : e_i < e_{i+1}\}$ and $\ROW(e) = \{e_i : 1 \leq i \leq n\} - \{1\}$. Note that the strict inequality in the definition of $\ASC$ is essential, and that we always have $e_1 = 1$. We let $\des = \# \DES$ be the number of descents, and similarly for $\ides$, $\asc$ and $\row$. If $S$ is a set of postive integers, we let $u^S = \prod_{i\in S} u_i$. From now on, we use indeterminates $u_1, u_2, \dots$ and $v$. 

\subsection*{Examples}
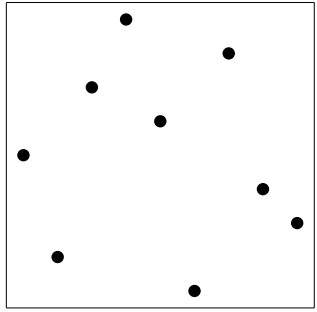
\begin{figure}
\scalebox{0.45} {
\begin{tikzpicture}
	\exI
\end{tikzpicture}
}
\caption{The child $\pi = 527961843$ referred to in Figure \ref{fi_pip}. We have $\ides(\pi) = \ides(\pi') + (2-1)_+ + 1 = 5$.}
\label{fi_pi}
\end{figure}

	We give the values of the statistics of the permutation $\pi$ in Figure \ref{fi_pi} and of the inversion table $e$ in Figure \ref{fi_e}. We have $\DES(\pi) = \{1,4,5,7,8\}$, $\IDES(\pi) = 1,3,4,6,8\}$, $\ASC(e) = \{1,2,4,5,7,8\}$ and $\ROW(e) = \{2,3,4,5,7\}$. Thus $\des(\pi) = 5$, $\ides(\pi) = 5$, $\asc(e) = 6$ and $\row(e) = 5$.

To prove Theorem \ref{th_main}, we will prove the stronger statement
\begin{equation}
	\label{eq_one}
	\sum_{\pi\in\mbb{S}_n} u^{\DES(\pi)}v^{\ides(\pi)} = \sum_{e\in\mbb{I}_n} u^{\ASC(e)}v^{\row(e)}.
\end{equation}

By M\"obius inversion, equation \ref{eq_one} holds if and only if we have 
\begin{equation}
	\label{eq_two}
	\sum_{\substack{\pi\in\mathbb{S}_n \\ S \subseteq \DES(\pi)}} s^{\ides(\pi)} = 
	\sum_{\substack{e\in\mbb{I}_n \\ S \subseteq \ASC(e)}} s^{\row(e)},
\end{equation}
	for each $S \subseteq [n-1]$. The idea now is to fix a subset $S$ of postive integers and induct on $n$ (in a sense to be specified) to prove equation \ref{eq_two}.

Thus fix a subset $S \subseteq \{1, 2, \dots\}$. We define two rooted labeled trees $T_{\mbb{S}}$ and $T_{\mbb{I}}$, as follows. The vertices of $T_{\mbb{S}}$ are all permutations $w_1\dots w_n$ whose length $n$ satisfies $n \notin S$ (in particular, the empty permutation, of length $0$, is a node). A permutation $\pi$ of length $r+s$ is a child of another permutation $\pi'$ of length $r$, where $s$ is smallest such that $s \geq 1$ and $r+s \notin S$, if the first $r$ letters of $\pi$ induce\footnote{that is, if we renumber the first $r$ letters of $\pi$ by $1,\dots,r$, we get the permutation $\pi'$} the same permutation as $\pi'$. It follows that the empty permutation is the root node. Each node in $T_{\mbb{S}}$ is labeled by the pair $(n, k)$, where $n$ is the length of the permutation and $k$ is its number of inverse descents.

Similarly, the vertices of $T_{\mbb{I}}$ are all inversion tables $e_1\dots e_n$ whose length $n$ satisfies $n\notin S$, and an inversion table $e$ of length $r+s$ is a child of $e'=e'_1\dots e'_r$ if $e = e'_1\dots e'_r e_{r+1} \dots e_{r+s}$ (with the same condition on $s$ as for $T_{\mathbb{S}}$). Each node is labeled $(n, k)$, where $n$ is the length of the inversion table and $k$ is its value of the row statistic. The empty inversion table is the root node.

To prove equation \ref{eq_two} for our fixed set $S$, it suffices to prove that $T_{\mbb{S}}$ and $T_{\mbb{I}}$ are isomorphic as labeled rooted trees. We will do this by producing an isomorphism $\Phi$ between the two trees, which takes a permutation $\pi$ of length $n$ with $k$ inverse descents satisfying $S \cap [n] \subseteq \DES(\pi)$ to some inversion table $e = \Phi(\pi)$ of length $n$, satisfying $\row(e) = k$ and $S \cap [n] \subseteq \ASC(e)$.

We will construct $\Phi$ inductively. Let $\Phi$ map the root of $T_{\mbb{S}}$ to the root of $T_{\mbb{I}}$. Suppose that we have already defined $\Phi(\pi') = e'$. We will show, for each $k$, that the number of children of $\pi'$ with $k$ inverse descents equals the number of children of $e'$ whose row statistic equals $k$. This allows us to extend $\Phi$ to all the children of $\pi'$. Thus fix a permutation $\pi'$ of length $r$ and an inversion table $e'$ of length $r$ such that $\ides(\pi') = \row(e') = p$, say. Suppose $s$ is smallest such that $s \geq 1$ and $r+s \notin S$. For children $\pi$ of $\pi'$, we call the first $r$ letters the {\it early} part, and the last $s$ letters the {\it late} part. 
The children $\pi$ are determined in a bijective way by $(r+1)$-tuples $(x_0, \dots, x_r)$ of nonnegative integers $x_i$ with sum $s$. The bijection\footnote{To prove this is a bijection, note that the last $r$ letters of $\pi$ form a decreasing word since $\pi$ is a child (which implies that $r, r+1, \dots, r+s-1$ are descents).} is given by letting $x_i$ be the number of letters in the late part of $\pi$ which are between (in value) the $i$th and $(i+1)$st largest letters of the early part of $\pi$. Moreover, the number of inverse descents of $\pi$ is\footnote{We use the notation $x_+ = x$ for $x \geq 0$ and $x_+ = 0$ for $x < 0$.} $\ides(\pi') + \sum_{i \in T} (x_i - 1)_+ + \sum_{i\notin T} x_i$, where $T = \IDES(\pi')$. This is hopefully made clear by Figures \ref{fi_pip} and \ref{fi_pi}.

On the other hand we consider children $e$ of $e'$, of length $r+s$. Such $e$ are in bijection with subsets $T \subseteq [r+s]$ of size $s$ by letting $T = \{e_i : r < i \leq r+s\}$ (note that this is a {\it set} since $\{s+1, \dots, r-1\} \subseteq \ASC(e)$ is a child). Moreover, $\row(e) = \row(e') + \#(T^c \cap ROW(e'))$.

By the preceding two paragraphs, all that remains is to prove that the number of $(x_0,\dots,x_r)$ with $\sum_{i=0}^{p-1} (x_i-1)_+ + \sum_{i = p} ^r x_i = t$ and $\sum_{i=0} ^r x_i = s$ (here, we have reordered the $x_i$'s, which clearly does not affect the count) equals the number of $T \subseteq [r+s]$ such that $|T| = s$ and $|T \cap \{p+2, p+3, \dots\}| = t$, for all nonnegative integers $r,s,p,t$.

These two counts are easily seen to be
\[
\sum_{a=0}^t {p+1\choose s-t} {s-a-1\choose a}{a+r-p-1\choose a} 
\]
and 
\[
{p+1\choose s-t}{r+s-p-1\choose t}
\] respectively. That they are equal is a classical fact. This finishes the proof of Theorem \ref{th_main}.

\subsection*{Final remarks}

It is interesting to note that while the statement of Theorem \ref{th_main} is symmetric in $t$ and $s$, the proof is not. We have failed to generalize the Theorem to one with two set-valued statistics. If we define $\maj(\pi)$ to be the sum of the elements in $\DES(\pi)$, and $\amaj(e)$ to be the sum of the elements in $\ASC(e)$, then the proof shows that the pairs $(\maj, \ides)$ and $(\amaj, \row)$ are equidistributed. There does not seem to be an obvious generalisation involving $\maj$, $\amaj$ and sums of $\IDES$ and $\ROW$. Finally, it follows from Theorem \ref{th_main} that $\asc$ and $\row$ are equidistributed (since $\des$ and $\ides$ are). I do not know of a direct proof of this fact.

\begin{thebibliography}{9}
\bibitem{BR} Matthias Beck and Sinai Robins: Computing the continuous discretely: Integer-point enumeration in polyhedra, Undergraduate Texts in Mathematics, Springer, New York, 2007.
\bibitem{B} Petter Br\"{a}nd\'{e}n: Actions on permutations and unimodality of descent polynomials. European J. Combin., 29(2):514-531, 2008.
\bibitem{P} Kyle Petersen: Two-sided Eulerian numbers via balls in boxes. Math. Mag. (to appear). Preprint available at arXiv:1209.6273v1, 2012.
\bibitem{S} Richard Stanley: Enumerative Combinatorics, Vol 2, Cambridge, 1999.
\bibitem{V} Mirk\'{o} Visontai: Some remarks on the joint distribution of descents and inverse descents, Electron. J. Combin., Volume 20, Issue 1, 2013, Research article 52, 12pp.
\end{thebibliography}
\end{document}